\documentclass[12pt]{article}
\usepackage{a4,amsmath,amsthm}
\usepackage{eucal}
\usepackage{pb-diagram}
\DeclareMathAlphabet{\bb}{U}{msb}{m}{n}
\gdef\C{\bb C}

\gdef\dR{\bb R}

\DeclareMathOperator{\spin}{{\bf Spin}}

\DeclareMathOperator{\Ker}{Ker}

\DeclareMathOperator{\sech}{{\rm sech}}

\newcommand{\re}{\mbox{\rm Re}\,}

\newcommand{\cA}{\mathcal{A}}

\newcommand{\M}{{\bf\sf M}}

\newcommand{\cl}{C\kern -0.2em \ell}

\newcommand{\p}{\prime}
\newcommand{\e}{\mbox{\bf e}}
\newcommand{\BH}{{\rm\hskip 0.1pt %
                            I\hskip -2.15pt H}} 
\newcommand{\R}{{\rm\hskip 0.1pt %
                           I\hskip -2.15pt R}}

\begin{document}
\title{Spinor Fields on the Surface of Revolution and their Integrable
Deformations via the mKdV-Hierarchy}
\author{Vadim V. Varlamov\\
{\it\small Siberia State University of Industry, 654007 Novokuznetsk, Russia}}
\date{}
\maketitle
\begin{abstract}
Spinor fields on surfaces of revolution conformally immersed into 3-dimensional
space are considered in the framework of the spinor representations of
surfaces. It is shown that a linear problem (a 2-dimensional Dirac equation)
related with a modified Veselov-Novikov hierarchy in the case of the surface
of revolution reduces to a well-known Zakharov-Shabat system. In the case
of one-soliton solution an explicit form of the spinor fields is given
by means of linear Bargmann potentials and is expressed via the Jost
functions of the Zakharov-Shabat system. It is shown also that integrable
deformations of the spinor fields on the surface of revolution are defined
by a modified Korteweg-de Vries hierarchy.
\end{abstract}
\section{Introduction}
At present time methods of the theory of surfaces penetrate into many
areas of theoretical and mathematical physics and become an important
and inherent part of the modern physics. The deep relation between the
theory of surfaces and soliton theory is well known, there exists a
numerous literature devoted to this question (see, for example, \cite{Sym85,
Bob94}). Recently, Konopelchenko and Taimanov introduced a new approach
related the methods of integrable systems with a theory of surfaces
conformally immersed into 3- and 4-dimensional spaces \cite{Kon96a,KT95,KT96,
Tai97a,Tai97b}. The main tool of this approach is a generalized Weierstrass
representation for a conformal immersion of surfaces into $\R^3$ or 
$\R^4$, and also a linear problem related with this representation.
A consideration of the linear problem along with the Weierstrass 
representation allows to express integrable deformations of surfaces via
such hierarchies of nonlinear differential equations as a modified
Veselov-Novikov hierarchy, Davy-Stewartson hierarchy and so on. At present
this approach intensively used for the study of constant mean curvature
surfaces, Willmore surfaces, surfaces of revolution and in many other
problems related with differential geometry and physics \cite{Kon96a,Kon98,
KL98a,KL98b,KT95,KT96,KL97,Mat97a,Mat97b,Mat98,Tai97a,Tai97b,Tai97c,Tai98,
Var98a,Var98b,Yam99}. Further,
it is known \cite{Ken79,HO83,HO85,KL98b} that the Weierstrass representation is
closely related with a Gauss map which widely used in string theory
\cite{BP88,Bud86,BDF86,BR87,VP94}, 
since a string worldsheet may be considered as a surface
conformally immersed into $\R^4$. Thus, we have a very intriguing relation
between string theory and soliton theory \cite{CK96,KL98c,KL99}.

On the other hand, there exists a natural formulation of the Weierstrass
representation in terms of spinor representation of surfaces 
\cite{At71,John80,Sul89,
Abr89,KS96}. By virtue of this we have an interesting relation with spinor
structures and spinor bundles, and also with a theory of the Dirac operator
on spin manifolds \cite{Bau81,LM89,Bar91,BFGK,Fr97a,Amm98}. In this direction,
spinor representations of surfaces immersed into 3-dimensional spin
manifold, and also upper bounds for the first eigenvalue of the Dirac
operator on such surfaces have been considered in the recent papers
\cite{Fr97b,AF98,AAF98}.

In the present paper we consider a relation between a Weierstrass
representation in the Eisenhardt-Konopelchenko form \cite{Eisen,Kon96a} for
a conformal immersion of surfaces into a 3-dimensional euclidean space
and spinor representations of surfaces. The linear problem related with
the mVN-hierarchy in this case is coincide with a 2-dimensional Dirac
operator. The spinor field $(\psi,\varphi)$ satisfying the 2-dimensional
Dirac equation is understood as a smooth section of the spinor bundle
over the Riemann surface. Integrable deformations of $(\psi,\varphi)$
are defined by the mVN-hierarchy. In parallel with this we consider
spinor representations of surfaces in terms of the theory of spin manifolds
and establish an algebraic structure of Weyl half-spinors $\varphi^+,
\varphi^-$ of the 2-dimensional manifold (surface) immersed into
a 3-dimensional spin manifold. It is shown that integrable deformations
of the field $(\varphi^+,\varphi^-)$ are also expressed via the
mVN-hierarchy. Further, in section 4 we consider a surface of revolution
and show that in this case the linear problem reduces to a well-known
Zakharov-Shabat system \cite{ZS71}, the fundamental solutions (Jost
functions) of which in the case of one-soliton solution are expressed
via the linear Bargmann potentials \cite{Bar49}. An explicit form of the
spinor fields $(\psi,\varphi),\,(\varphi^+,\varphi^-)$ is given.
In conclusion, using a well-known fact that in the case of the surface of
revolution the mVN-hierarchy reduces to a modified Korteweg-de Vries
(mKdV) hierarchy \cite{KT95,Tai97a} we give an explicit dependence of the
spinor fields on deformation variable $t$ in the case of the one-soliton
solution and the first member of the mKdV-hierarchy.
\section{Spinor structures on manifolds and the Dirac operator on surfaces}
Let us recall some basic facts about spinor structures on manifolds.
Let $(M^{p,q},g)$ be $n$-dimensional $(n=p+q)$ oriented (time-oriented and
space-oriented) pseudo-Riemann manifold and let $P\times_{SO(p,q)}$
be the bundle
of all $SO(p,q)$-frames of $(M^{p,q},g)$. Then, {\it a spinor structure}
of $(M^{p,q},g)$ is a pair $(Q,f)$ of a $\spin(p,q)$-principal bundle
$Q\times_{\spin(p,q)}$ and a continuous surjective map $f:\,Q\longrightarrow
P$ such that the diagram
\[
\dgARROWLENGTH=0.75em
\begin{diagram}
\node{Q\times_{\spin(p,q)}\times\spin(p,q)}\arrow[2]{e}\arrow[2]{s,r}
{f\times\lambda}\node[2]{Q\times_{\spin(p,q)}}\arrow{se}\arrow[2]{s,r}{f}
\\
\node[4]{M}\\
\node{P\times_{SO(p,q)}\times SO(p,q)}\arrow[2]{e}\node[2]{P\times_{SO(p,q)}}
\arrow{ne}
\end{diagram}
\]
commutes, here $\lambda$ is a double covering $\spin(p,q)\longrightarrow
SO(p,q)$.

Further, an oriented pseudo-Riemann manifold with a spinor structure
$(Q,f)$ is called {\it a spin manifold}.
The complex vector bundle
\[
S=Q\times_{\spin(p,q)}\Delta_{p,q}
\]
associated with the $\spin(p,q)$-principal bundle $Q\times_{\spin(p,q)}$
by means of the spinor representation $\Delta_{p,q}$ is called {\it
a spinor bundle} of the manifold $(M^{p,q},g)$. A smooth section
$\varphi\in\Gamma(S)$ of the bundle $S$ is called {\it a spinor field}
on the manifold $(M^{p,q},g)$. In case of even dimension, $p-q\equiv
0,2,4,6\pmod{8}$, the spinor bundle $S$ splits into two subbundles
$S=S^+\oplus S^-$, where
\[
S^{\pm}=Q\times_{\spin(p,q)}\Delta^{\pm}_{p,q},
\]
here $\Delta^{\pm}_{p,q}$ are the spaces of half-spinors (Weyl spinors).

Let us consider now a 3-dimensional oriented Riemann manifold $M^{3,0}$
with a fixed spin structure and also a surface $M^{2,0}$ conformally
immersed into $M^{3,0}$. The surface $M^{2,0}$ is considered as a 
space-oriented Riemann manifold with a spin structure. Moreover, since
the normal bundle of $M^{2,0}$ is trivial, then the spin structure of
$M^{3,0}$ induces a spin structure on the surface $M^{2,0}$. Let $\Phi$
be a spinor field on the manifold $M^{3,0}$. It is obvious that the field
$\Phi$ is a section of 2-dimensional spinor bundle $S=Q\times_{\spin(3,0)}
\Delta_{3,0}$, since $\spin(3,0)\cong\spin(0,3)$ and $\cl^+_{0,3}\cong\BH$,
then $\dim\Delta_{3,0}=2$ (here $\BH$ is a quaternion algebra, $\cl^+_{p,q}$
is a Clifford algebra of all even elements). On the other hand,
$\spin(3,0)\cong Sp(1)\cong SU(2)\cong S^3$ \cite{Port} and $\cl_{3,0}
\cong\C_2\cong\M_2(\C)$. Further, it takes to find a restriction of the
spinor bundle $S=Q\times_{\spin(3,0)}\Delta_{3,0}$ of the manifold $M^{3,0}$
onto a spinor bundle $S_{M^{2,0}}=Q\times_{\spin(2,0)}\Delta_{2,0}$ of the
surface $M^{2,0}$ conformally immersed into $M^{3,0}$. Let $\phi$ is a
spinor field on the surface $M^{2,0}$. Obviously, this field is a section
of 2-dimensional spinor bundle $S=Q\times_{\spin(2,0)}\Delta_{2,0}$, since
$\spin(2,0)\cong\spin(0,2)$ and $\cl^+_{0,2}\cong\cl_{0,1}\cong\C$, then
$\spin(2,0)\cong U(1)\cong S^1$. Moreover, the spinor bundle of the surface
$M^{2,0}$ splits into two subbundles,
\[
S=S^+\oplus S^-,
\]
where $S^{\pm}=Q\times_{\spin(2,0)}\Delta^{\pm}_{2,0}$. Respectively,
a smooth section $\phi\in\Gamma(S)$ of the bundle $S$ has a form
$\phi=\phi^++\phi^-$, where \cite{Fr97b}
\begin{equation}\label{e1}
\phi^+=\frac{1}{2}(\phi+i\vec{N}\phi),\quad\phi^-=\frac{1}{2}(\phi-
i\vec{N}\phi),
\end{equation}
here $\vec{N}=\e_3=\e_1\e_2$.\\[0.3cm]
{\bf Remark}. Let us consider in more details algebraic definition of the
spinors (\ref{e1}). At first, for the 3-dimensional Clifford algebra
$\cl_{3,0}$ associated with the manifold $M^{3,0}$ there exists an
isomorphism $\cl_{3,0}\cong\C_2$, where $\C_2$ is an algebra of complex
quaternions. Further, according to \cite{Che54,Lou81} spinors are
considered as elements of the minimal left ideal $I_{p,q}=\cl_{p,q}e_{pq}$,
where $e_{pq}$ is a primitive idempotent $e_{pq}=\frac{1}{2}(1+\e_{\alpha_1})
\ldots\frac{1}{2}(1+\e_{\alpha_k}),\,k=q-r_{q-p},\,r_{q-p}$ are the
Radon-Hurwitz numbers defined by the recurrence formula $r_{i+8}=r_i+4$ and
\begin{center}
\begin{tabular}{lcccccccc}
$i$   & 0 & 1 & 2 & 3 & 4 & 5 & 6 & 7\\ \hline
$r_i$ & 0 & 1 & 2 & 2 & 3 & 3 & 3 & 3
\end{tabular}
\end{center}
The primitive idempotent of $\cl_{3,0}\cong\C_2$ has a form $e_{30}=
\frac{1}{2}(1+\e_0)\sim\frac{1}{2}(1+i\e_{12})$, since in this case
$k=q-r_{q-p}=0-r_{-2}=0-(r_6-4)=1$. It is obvious that $I_{3,0}=\cl_{3,0}e_{30}
\cong\C_2e_{30}\cong\M_2(\C)e_{30}$, therefore for the algebraic definition
of the spinors (\ref{e1}) we have
\begin{eqnarray}
\varphi^+&=&e^+_{30}I_{3,0},\nonumber\\
\varphi^-&=&e^-_{30}I_{3,0},\nonumber
\end{eqnarray}
here
\[
e^+_{30}=\frac{1}{2}(1+i\e_{12}),\quad e^-_{30}=\frac{1}{2}(1-i\e_{12})
\]
are the mutually orthogonal idempotents.\\[0.3cm]
Further, according to \cite{Fr97b} suppose that the spinor $\Phi$ on the
manifold $M^{3,0}$ is a real Killing spinor, i.e. there exists a number
$\lambda\in\R$ such that for any tangent vector $\vec{T}\in\tau_{M^{3,0}}$
the derivative of the spinor $\Phi$ in the direction of $\vec{T}$ is
given by the equation
\[
\bigtriangledown^{M^{3,0}}_{\vec{T}}\Phi=\lambda\cdot\vec{T}\cdot\Phi.
\]
Then the Dirac equation for the spinors (\ref{e1}) may be written as
follows \cite{Fr97b}
\begin{eqnarray}
D\varphi^-&=&(-2\lambda+iH)\varphi^+,\nonumber\\
D\varphi^+&=&(-2\lambda-iH)\varphi^-,\nonumber
\end{eqnarray}
where $H$ is a mean curvature of the surface $M^{2,0}$ and
\[
D\varphi^\pm=\e_1\cdot\bigtriangledown^{M^{2,0}}_{\e_1}\varphi^\pm+
\e_2\cdot\bigtriangledown^{M^{2,0}}_{\e_2}\varphi^\pm
\]
is a Dirac operator of $M^{2,0}$. On the other hand, if $\Phi$ is a
parallel spinor on the manifold $M^{3,0}$ ($\lambda=0$), then
\begin{equation}\label{e2}
\begin{array}{ccc}
D\varphi^-&=&iH\varphi^+,\\
D\varphi^+&=&-iH\varphi^-.
\end{array}
\end{equation}
%Finally, for the spinor field
%\begin{equation}\label{e3}
%\varphi^\ast=\varphi^+-i\varphi^-=\frac{1}{2}(1-i)\Phi_{|M^{2,0}}+
%\frac{1}{2}(-1+i)\vec{N}\cdot\Phi_{|M^{2,0}}
%\end{equation}
%of constant length on $M^{2,0}$ satisfying the two-dimensional Dirac
%equation
%\[
%D\varphi^\ast=H\varphi^\ast
%\]
%the following theorem take place
%\begin{theorem} [{\rm Agricola-Friedrich}]
%The first eigenvalue $\lambda^2_1$ of the square of the Dirac operator
%on a surface $M^{2,0}\hookrightarrow M^{3,0}$ is bounded by
%\[
%\lambda^2_1\leq\frac{\int_{M^{2,0}}H^2(f^2+G^2(f))dM^2+\int_{M^{2,0}}
%|grad\,f|^2(1+[G^{\prime}(f)^2]dM^2}{\int_{M^{2,0}}(f^2+G^2(f))dM^2}
%\]
%where $f:\,M^{2,0}\rightarrow\R$, $G:\,\R\rightarrow\R$ are smooth
%functions and $G^\prime$ denotes the derivative of $G$.
%\end{theorem}
\section{The spinor and Weierstrass representations of surfaces and
the modified Veselov-Novikov hierarchy} 
Let $S$ be a connected Riemann surface with a local complex coordinate $z$
and let $P\times_{G}$ be a principal bundle of $S$ with the structure
group $G$ ($G$ is a group of fractional linear transformations). Then the
{\it spinor representation of a surface} is given by the following
diagram
\[
\dgARROWLENGTH=0.8em
\begin{diagram}
\node{Q\times_{\widetilde{G}}}\arrow{s,l}{\mu}\arrow{e,t}{\chi}
\node{Q\times_{\spin(2,0)}}\arrow{s,r}{f}\\
\node{P\times_{G}}\arrow{s}\arrow{e,t}{\omega}\node{P\times_{SO(2,0)}}
\arrow{s}\\
\node{S}\arrow{e,t}{g}\node{M^{2,0}}
\end{diagram}
\]
Here $g:\,S\longrightarrow\C\cup\{\infty\}$ is a Gauss map, which can be
viewed as an identification (via the stereographic projection) of the
surface $M^{2,0}$ and $\C\cup\{\infty\}$, $\omega$ is a bundle map of
$P\times_{G}$ into $P\times_{SO(2,0)}$ satisfying the integrability
condition $\re d\omega=0$ \cite{KS96}. $\widetilde{G}$ is a double covering
of $G$, i.e. there is a mapping $\lambda:\,\widetilde{G}\longrightarrow G$.
According to \cite{Roz55} the group $\widetilde{G}$ is a group of complex
matrices of the form $\begin{pmatrix} A & B \\ -\bar{B} & \bar{A}
\end{pmatrix}$, and $\Ker\lambda$ consists of the two matrices
$\begin{pmatrix} 1 & 0 \\ 0 & 1\end{pmatrix}$ and $\begin{pmatrix}
-1 & 0\\ 0 & -1 \end{pmatrix}$. Further, if $(\omega,g)$ is a bundle map
of $P\times_{G}$ into $P\times_{SO(2,0)}$, then $(i)$ there is a unique
structure $Q\times_{\widetilde{G}}$ on $S$ such that $\omega$ lifts to a
bundle map $\chi:\,Q\times_{\widetilde{G}}\longrightarrow Q\times_{\spin(2,0)}$;
$(ii)$ there are exactly two such lifts $\psi$, and these differ only
by sign.

On the other hand, the Weierstrass representation \cite{Weier} describes
a conformal immersion of minimal surfaces into 3-dimensional Euclidean space:
\[
f=\re\left[\int(1-s^2,i(1+s^2),2s)\mu(s)\right]:\,M^{2,0}\longrightarrow
\R^{3,0},
\]
where
\[
s=\frac{H(w)}{G(w)},\quad\mu(s)=G^2\frac{dW}{ds},
\]
and $H(w),G(w)$ are holomorphic functions defined in a circle or in all
complex plane. The Weierstrass and spinor representations are related by
the equation \cite{KS96}
\[
f=\sigma(\psi)=(\varphi^2_1-\varphi^2_2,i(\varphi^2_1+\varphi^2_2),
2\varphi_1\varphi_2),
\]
where $\psi=(\varphi_1,\varphi_2)$ is a pair of sections of the bundle
$Q\times_{\widetilde{G}}$, and
\begin{equation}\label{e4}
\mu(s)=\varphi^2_1,\quad s=\varphi_2/\varphi_1.
\end{equation}
In connection with this let us consider a Weierstrass representation
in the Eisenhardt-Konopelchenko form \cite{Eisen,Kon96a}:
\begin{eqnarray}
X^1&+&iX^2=i\int_\Gamma(\bar{\psi}^2dz^{\p}-\bar{\varphi}^2d\bar{z}^{\p}),
\nonumber\\
X^1&-&iX^2=i\int_\Gamma(\varphi^2dz^{\p}-\psi^2d\bar{z}^{\p}),\nonumber\\
X^3&=&-\int_\Gamma(\psi\bar{\varphi}dz^{\p}+\varphi\bar{\psi}d\bar{z}^{\p}),
\label{e5}
\end{eqnarray}
where $\Gamma$ is an arbitrary curve in the complex plane $\C$, $\psi$
and $\phi$ are complex-valued functions on variables $z,\bar{z}\in\C$
satisfying the following linear system (2-dimensional Dirac equation):
\begin{equation}\label{e6}
\begin{array}{ccc}
\psi_z&=&p\varphi,\\
\varphi_{\bar{z}}&=&-p\psi,
\end{array}
\end{equation}
where $p(z,\bar{z})$ is a real-valued function, and $p=\frac{u}{2}H$,
where $H$ is a mean curvature of the surface, $u=|\psi|^2+|\varphi|^2$.
If to interpret the functions $X^i(z,\bar{z})$ as coordinates in the
space $\R^{3,0}$, then the formulae (\ref{e5}),(\ref{e6}) are defined
a conformal immersion of the surface $M^{2,0}$ into $\R^{3,0}$.
The formulae (\ref{e5}) reduce to a classical Weierstrass representation
if suppose
\[
\varphi={\rm G},\quad\psi={\rm H}. 
\]

Further, it is well-known that the mVN-hierarchy is related with the
system (\ref{e6}), see \cite{Kon96a,KT95,KT96,Tai97a}. This hierarchy is
appear as compatibility conditions
\begin{equation}\label{e7}
\left[\frac{\partial}{\partial t_n}-A_n,L\right]=B_nL,
\end{equation}
where the Lax operator $L$ and deformation operators $A_n,B_n$ are defined
by the following expressions
\begin{gather}
L=\begin{pmatrix}
\partial & -p\\
p & \bar{\partial}
\end{pmatrix}\nonumber\\
A_n=\partial^{2n+1}+\sum^{2n-1}_{i=0}X^i\partial^i+\bar{\partial}^{2n+1}
+\sum^{2n-1}_{i=0}\widetilde{X}^i\bar{\partial}^i,\nonumber\\
B_n=\sum^{2n-1}_{i=0}S^i\partial^i+\sum^{2n-1}_{i=0}\widetilde{S}^i
\bar{\partial}^i,\nonumber
\end{gather}
here $X^i,\widetilde{X}^j,S^i,\widetilde{S}^j\;(i,j=0,\ldots,2n-1)$ are
$2\times 2$ deformation matrices. The compatibility conditions (\ref{e7})
give the following deformation equations
\[
\frac{\partial}{\partial t_n}p=\partial^{2n+1}p+\bar{\partial}^{2n+1}p+
\ldots,
\]
which consist the mVN-hierarchy. The deformation of the eigenfunctions
of the operator $L$ are defined by the system
\begin{equation}\label{e8}
\frac{\partial}{\partial t_n}\begin{pmatrix} \psi \\ \varphi\end{pmatrix}=
A_n\begin{pmatrix} \psi \\ \varphi\end{pmatrix}.
\end{equation}
When $n=1$ we have the first member of the
mVN-hierarchy:
\begin{equation}\label{e7'}
p_t=p_{t^+}+p_{t^-},
\end{equation}
where
\begin{eqnarray}
p_{t^+}&=&p_{zzz}+3p_z\omega+\frac{3}{2}p\omega_z,\nonumber\\
p_{t^-}&=&p_{\bar{z}\bar{z}\bar{z}}+3p_{\bar{z}}\bar{\omega}+\frac{3}{2}p
\bar{\omega}_{\bar{z}}.\nonumber
\end{eqnarray}
Further, the deformation of the eigenfunctions $\psi,\varphi$ via (\ref{e8})
generates the corresponding deformations of the coordinates $X(z,\bar{z},t)$
of the surface conformally immersed into space $\R^{3,0}$.

It is easy to verify that the matrices of the group $\widetilde{G}$, which
double covers the group of the fractional linear transformations of the
Riemann surface $S$, are coincide with matrices of a matrix representation
of the quaternion algebra $\cl_{0,2}=\BH$. Indeed, a general element of
$\BH$ has a form
\[
\cA=a^0\e_0+a^1\e_1+a^2\e_2+a^{12}\e_{12},
\]
where $a^0,a^1,a^2,a^{12}\in\R$. Further, in virtue of the maps
$\e_i\rightarrow\sigma_i\;(i=0,1,2)$, where
\[
\sigma_0=\begin{pmatrix}
1 & 0\\
0 & 1
\end{pmatrix},\quad\sigma_1=\begin{pmatrix}
0 & 1\\
-1 & 0
\end{pmatrix},\quad\sigma_2=\begin{pmatrix}
0 & i\\
i & 0
\end{pmatrix},
\]
we have
\[
a^0\sigma_0+a^1\sigma_1+a^2\sigma_2+a^{12}\sigma_{12}=\begin{pmatrix}
A & B\\
-\bar{B} & \bar{A}
\end{pmatrix},
\]
where 
\[
A=a^1+ia^{12},\quad B=a^1+ia^2.
\]
Therefore, $\spin(2)\subset\widetilde{G}$. 
Let $Q\times_{\widetilde{G}}\Delta_2$
be {\it a spinor bundle of the Riemann surface} $S$, where $\Delta_2$ is a
spinor representation of the group $\widetilde{G}$, and let $(\psi,\varphi)$
are smooth sections of the bundle $Q\times_{\widetilde{G}}\Delta_2$,
where the spinors $\psi,\varphi$ satisfy the Dirac equation (\ref{e6}). 
Further, let
$\beta:\,Q\times_{\widetilde{G}}\Delta_2\longrightarrow Q\times_{\spin(2,0)}
\Delta_{2,0}$ be a spinor bundle map between the spinor bundle of $S$ and
spinor bundle of the manifold $M^{2,0}$. Then in virtue of the map $\beta$
we can to identify the sections (spinor fields) of these bundles
\[
\begin{pmatrix}
\psi\\
\varphi
\end{pmatrix}\sim\begin{pmatrix}
\varphi^+\\
\varphi^-
\end{pmatrix}.
\]
In so doing, $\e_1\rightarrow 1\cdot x,\,\e_2\rightarrow i\cdot y$, and
$\e_1\cdot\bigtriangledown^{M^{2,0}}_{\e_1}\rightarrow\frac{\partial}
{\partial x},\,\e_2\cdot\bigtriangledown^{M^{2,0}}_{\e_2}\rightarrow
i\frac{\partial}{\partial y}$. Therefore, for the Dirac operator of
$M^{2,0}$ we have
\[
D\sim\frac{1}{2}\left(\frac{\partial}{\partial x}+i\frac{\partial}
{\partial y}\right).
\]
Thus, in virtue of the spinor bundle map $\beta:\,Q\times_{\widetilde{G}}
\Delta_2\longrightarrow Q\times_{\spin(2,0)}\Delta_{2,0}$ the systems
(\ref{e2}) and (\ref{e6}) are equivalent. Therefore, {\it if $\Phi$ is a
parallel spinor field on the manifold $M^{3,0}$, then integrable
deformations of the restriction $\Phi_{|M^{2,0}}=\begin{pmatrix}
\varphi^+\\ \varphi^-\end{pmatrix}$, where $\varphi^+$ and $\varphi^-$
are the half-spinors (Weyl spinors) of the surface $M^{2,0}$, 
are defined by the
mVN-hierarchy}. 
\section{Surfaces of revolution}
Let us consider now one important particular case, namely, the case when
the surface $M^{2,0}$ is a surface of revolution. In this case
\cite{KT95,KT96,Tai97a,Tai97b} the functions $\psi$ and $\varphi$ take
the form
\begin{equation}\label{e9}
\psi=r(x)\exp(\lambda y),\quad\varphi=s(x)\exp(\lambda y),
\end{equation}
and the system (\ref{e6})
\begin{eqnarray}
\frac{\partial}{\partial z}\psi&=&p\varphi,\nonumber\\
\frac{\partial}{\partial\bar{z}}\varphi&=&-p\psi,\nonumber
\end{eqnarray}
where
\[
\frac{\partial}{\partial z}=\frac{1}{2}\left(\frac{\partial}{\partial x}+
i\frac{\partial}{\partial y}\right),\quad\frac{\partial}{\partial\bar{z}}=
\frac{1}{2}\left(\frac{\partial}{\partial x}-i\frac{\partial}{\partial y}
\right),
\]
reduces to a system
\begin{eqnarray}
\frac{\partial}{\partial x}r+i\lambda r&=&us,\nonumber\\
\frac{\partial}{\partial x}s-i\lambda s&=&-ur,\label{e10}
\end{eqnarray}
where $u=2p$, and $p(x)$ is a real-valued function. The system (\ref{e10})
is nothing but a well-known {\it Zakharov-Shabat system} \cite{ZS71}.
Solutions of this system can be obtained via the Bargmann potentials
\cite{Bar49} (see also \cite{Lam80}):
\begin{eqnarray}
r&=&e^{-i\lambda x}[4i\lambda+a(x)],\nonumber\\
s&=&e^{-i\lambda x}b(x).\label{e11}
\end{eqnarray}
Indeed, substituting the Bargmann potentials (\ref{e11}) into (\ref{e10})
we find
\begin{equation}\label{e12}
a_x=ub,\quad b=2u,\quad b_x=-ua.
\end{equation}
Excluding $u$ from the first and third equations (\ref{e12}) and then
integrating it we obtain $a^2+b^2=4\mu^2$, where $\mu$ is a constant of
integration. Further, from the first and second equations (\ref{e12}) we
have $b^2=2a_x$. Thus, the function $a(x)$ satisfies the equation
\[
a_x+a^2/2=2\mu^2.
\]
After substitution $a=2w_x/w$ we came to the linear equation
\begin{equation}\label{e13}
w_{xx}-\mu^2w=0.
\end{equation}
The solution of the equation (\ref{e13}) has a form
\[
w=\alpha e^{\mu x}+\beta e^{-\mu x}=2e^{\phi}\cosh(\mu-\phi),
\]
where $\phi=\frac{1}{2}\ln(\beta/\alpha)$. Then from the first and second
equations (\ref{e12}) we obtain
\[
u^2=\frac{1}{2}a_x=\mu^2\sech^2(\mu x-\phi)
\]
and finally
\begin{equation}\label{e13'}
u=\pm\mu\sech(\mu x-\phi).
\end{equation}
Let $\mu>0$, from (\ref{e12}) we find the functions $b(x)$ and $a(x)$:
\begin{gather}
b(x)=\pm 2\mu\sech(\mu-\phi),\nonumber\\
a(x)=4\mu^2\int\sech^2(\mu x-\phi)dx=4\mu^2\tanh(\mu x-\phi).\nonumber
\end{gather}
Substituting now the functions $a(x),b(x)$ into (\ref{e11}) we find
finally the fundamental solutions (Jost functions) of the system (\ref{e10}):
\begin{eqnarray}
\phi^+_1(x,\lambda)&=&e^{-i\lambda x}\frac{2i\lambda+\mu\tanh(\mu x-\phi)}
{2i\lambda-\mu},\nonumber\\
\phi^+_2(x,\lambda)&=&\pm e^{-i\lambda x}\frac{\sech(\mu x-\phi)}
{2i\lambda-\mu}.\label{e14}
\end{eqnarray}
Therefore, {\it in the case of surface of revolution the sections
$(\psi,\varphi)$ of the spinor bundle $Q\times_{\widetilde{G}}\Delta_2$
are expressed via the fundamental solutions of the Zakharov-Shabat
system as follows}
\begin{equation}\label{e15}
\begin{array}{ccc}
\psi&=&\phi^+_1(x,\lambda)e^{\lambda y},\\
\varphi&=&\phi^+_2(x,\lambda)e^{\lambda y}.
\end{array}
\end{equation}
Thus, in this case we obtain an explicit form of the spinors $\psi$ and
$\varphi$. Moreover, in virtue of the spinor bundle map 
$Q\times_{\widetilde{G}}\Delta_2\longrightarrow Q\times_{\spin(2,0)}
\Delta_{2,0}$ we can to express the half-spinors $\varphi^+,\varphi^-$
of the surface $M^{2,0}$ via the solutions (\ref{e14}).

Further, the mVN-hierarchy related with the system (\ref{e6}) in the
case of surface of revolution reduces to a modified Korteweg-de Vries
hierarchy. Indeed, for the first member (\ref{e7'}) of the mVN-hierarchy
the each parts $p_{t^+}$ and $p_{t^-}$ equals to $p_{xxx}+6p^2p_x$.
Take into account that $p=u/2$ we came to a famous modified
Korteweg-de Vries equation
\begin{equation}\label{e16}
u_t=u_{xxx}+\frac{3}{2}u^2u_x.
\end{equation}
Therefore, integrable deformations of the surfaces of revolution and also
a deformation of the spinor field $(\psi,\varphi)$ are defined by the
mKdV-hierarchy. Substitute now $u=\sqrt{4}v$ into (\ref{e16}) we obtain
a standard form of the mKdV-equation,
\begin{equation}\label{e17}
v_t-v_{xxx}-6v^2v_x=0.
\end{equation}
Further, plugging $v=u/\sqrt{4}$ into (\ref{e17}), where $u$ is a potential
of the form (\ref{e13'}), we see that a dependence on $t$ is expressed
by an equation
\[
\phi_t+\mu^3=0.
\]
Hence it immediately follows that $u=\pm\sech(\mu x-\mu^3t)$. Therefore,
an explicit dependence of $(\psi,\varphi)$ on deformation variable $t$
is defined by the expressions
\begin{eqnarray}
\psi&\sim&\phi^+_1(x,\lambda,t)e^{\lambda y}=e^{\lambda(y-ix)}
\frac{2i\lambda+\mu\tanh(\mu x-\mu^3t)}{2i\lambda-\mu},\nonumber\\
\varphi&\sim&\phi^+_2(x,\lambda,t)e^{\lambda y}=e^{\lambda(y-ix)}
\frac{\sech(\mu x-\mu^3t)}{2i\lambda-\mu}.\label{e18}
\end{eqnarray}
Thus, {\it integrable deformations of the spinor field $(\psi,\varphi)$
on the surface of revolution in the case of one-soliton solution of the
Zakharov-Shabat system (\ref{e10}) defining by the linear Bargmann
potentials (\ref{e11}), are expressed by the mKdV-equation and have 
an explicit form (\ref{e18})}.

Owing to absence of deformation, solutions (\ref{e15}) will be called
{\it static solutions}.

Further generalization is a transition to multi-soliton solutions
(for example, quadratic Bargmann potentials correspond to a two-soliton
solution), and also consideration of higher members of the mKdV-hierarchy,
which will be considered in a separate paper.

In conclusion, by virtue of the map $\beta:\,Q\times_{\widetilde{G}}
\longrightarrow Q\times_{\spin(2,0)}\Delta_{2,0}$ and relation between
smooth sections of these bundles, $(\psi,\varphi)\sim (\varphi^+,\varphi^-)$,
follows that {\it if $\Phi$ is a parallel spinor field on the manifold
$M^{3,0}$, and the surface $M^{2,0}$ conformally immersing into $M^{3,0}$
is homeomorphic to a surface of revolution, then integrable deformations
of the restriction $\Phi_{|M^{2,0}}=\begin{pmatrix}\varphi^+\\ \varphi^-
\end{pmatrix}$, where $\varphi^+$ and $\varphi^-$ are the Weyl spinors
of the surface $M^{2,0}$, are defined by the modified Korteweg-de Vries
hierarchy}.

\end{document}